\documentclass [12pt, a4paper] {article}

\usepackage {amsthm}
\theoremstyle {definition}
\newtheorem {definition} {Definition}

\newtheorem {theorem} [definition] {Theorem}
\newtheorem* {theorem*} {Theorem}
\newtheorem {proposition} [definition] {Proposition}
\newtheorem {corollary} [definition] {Corollary}

\theoremstyle {remark}

\usepackage {amssymb}

\usepackage {amsmath}

\DeclareMathOperator {\cl} {cl}

\renewcommand {\S} {\mathrm S}

\makeatletter
\newcommand {\forksym} {\raise0.2ex\hbox{\ooalign{\hidewidth$\vert$\hidewidth\cr\raise-0.9ex\hbox{$\smile$}}}}
\def\@forksym@#1#2{\mathrel{\mathop{\forksym}\displaylimits_{#2}}}
\def\nforks{\@ifnextchar_{\@forksym@}{\forksym}}

\makeatother

\usepackage {tikz}
\usetikzlibrary {calc}
\usetikzlibrary {intersections}
\usepackage {natbib}

\title {Groups with Locally Modular Homogeneous Pregeometries are Commutative}
\author {Levon Haykazyan}

\begin {document}
\maketitle

\begin {abstract}
It is well known that strongly minimal groups are commutative. Whether this is
true for various generalisations of strong minimality has been asked in several
different settings (see \cite{hyttinen, maesono, pilovic}). In this note we
show that the answer is positive for groups with locally modular homogeneous
pregeometries.
\end {abstract}

\section {Introduction}

The basic setting of geometric stability theory is a strongly minimal theory
where algebraic closure induces a pregeometry. In the literature there are
several approaches in generalising this basic setting. 

One approach is to ask for every definable set to be countable or cocountable
(instead of finite or cofinite in strongly minimal structures). An uncountable
structure satisfying this condition is called {\em quasiminimal}. 

Another option is to substitute the finite/cofinite dichotomy by large/small
dichotomy with respect to some ultrafilter on the boolean algebra of definable
sets (i.e. a type). So given a type $p \in \S(M)$ on a structure $M$, one can
look at the operator associating to $A \subseteq M$ the union of all small
$A$-definable sets (i.e. those that are not in $p$). The type $p$ is called {\em
strongly regular} if this operator is a closure operator (i.e. a pregeometry
possibly without exchange).

The third possibility is to simply require the existence of a homogeneous
pregeometry on a structure $M$. The pregeometry $\cl : \mathcal P(M) \to
\mathcal P(M)$ is called {\em homogeneous} if $M$ is infinite dimensional with
respect to $\cl$ and for every finite $A \subset M$ and $b, c \in M \setminus
\cl(A)$ there is an automorphism $f$ of $M$ (and of $\cl$) fixing $A$ pointwise
with $f(b) = c$.

By a well known result of \cite{reineke} all strongly minimal groups are
commutative. It has been asked whether the same is true in each of the above
three settings (in \cite{maesono} for quasiminimality, in \cite{pilovic} for
regular types and in \cite{hyttinen, hls} for groups with homogeneous
pregeometries). Reineke's argument has been generalised to show that if a
noncommutative group has one of the above properties, then it must have the
strict order property, so in particular is unstable. (See also \cite{gogacz}
which, as a step towards a potential counterexample, constructs an uncountable
group were all nontrivial elements are conjugate and have countable
centralisers.) However not much is known beyond this. The main result of this
paper is the following

\begin {theorem}
\label{main}
A group with a locally modular homogeneous pregeometry is commutative.
\end {theorem}

Using the connections between quasiminimality, strong regularity and homogeneous
pregeometries (see \cite{pilovic}), we have an analogues statement for
quasiminimal groups of cardinality at least $\aleph_2$ and groups with symmetric
strongly regular types, whose respective pregeometries are locally modular.

The initial idea of the proof was suggested to me by Alfonso Ruiz Guido. It used
the representation theorem of projective geometries and required a stronger
homogeneity hypothesis. Ivan Tomasic suggested to eliminate the use of the
representation theorem by finding the contradiction explicitly. Doing that
helped relax the homogeneity hypothesis. I express my gratitude to them and to
Rahim Moosa for encouraging me to write this up.

\section {Background on Pregeometries}

Although the motivation for Theorem \ref{main} comes from model theory, its
statement and proof involve none. In order to make the presentation accessible
to a wider audience, we provide a brief introduction to pregeometries.

\begin {definition}
A {\em pregeometry} (also called {\em matroid}) on a set $X$ is an operator $\cl
\colon \mathcal P(X) \to \mathcal P(X)$ satisfying the following conditions for
all $A \subseteq X$ and $a, b \in X$.
\begin {itemize}
\item Reflexivity: $A \subseteq \cl(A)$.
\item Transitivity: $\cl(\cl(A)) = \cl(A)$.
\item Finite character: $a \in \cl(A)$ if and only if $a \in \cl(A_0)$ for some
finite $A_0 \subseteq A$.
\item Exchange: $a \in \cl(A,b) \setminus \cl(A)$ implies $b \in \cl(A,a)$.
\end {itemize}
(Here and throughout we may use comma for union and omit the brackets for
singletons, so that e.g. $\cl(A,b)$ means $\cl(A \cup \{b\})$.)
\end {definition}

On any set $X$ the operator $\cl$ defined by $\cl(A) = A$ is a pregeometry.
Other examples of pregeometries include a vector space $V$ where $\cl(A)$ is the
linear span of $A$ and an algebraically closed field $F$ where $\cl(A)$ is the
algebraic closure of the field generated by $A$.

A subset $A \subseteq X$ is called {\em independent} if $a \not \in \cl(A
\setminus {a})$ for every $a \in A$. It is called {\em generating} if $\cl(A) =
X$ and a {\em basis} if it is both independent and generating. The exchange
axiom can be used to show that all bases have the same cardinality which is then
called the {\em dimension} of $X$ and denoted by $\dim(X)$.

Given $Y \subseteq X$ we can define the {\em restriction} $\cl^Y$ of $\cl$ onto
$Y$. This is a pregeometry on $Y$ defined by $\cl^Y(A) = \cl(A) \cap Y$. We also
define the {\em localisation} $\cl_Y$ of $\cl$ on $Y$. The localisation is a
pregeometry on $X$ defined by $\cl_Y(A) = \cl(A \cup Y)$. By $\dim(Y)$ we mean
the dimension of $\cl^Y$ and by $\dim(Y/Z)$ we mean the dimension of $Y$ in
$\cl_Z$.

A pregeometry $\cl$ is called a {\em geometry} if $\cl(\emptyset) = \emptyset$
and $\cl(a) = \{a\}$ for every $a \in A$. The idea is that closed subsets of
dimension $1$, $2$, $3$, etc correspond to points, lines, planes, etc. (Note
that what we call the dimension is $1$ more than the geometric dimension.) In
such a geometry one can only talk about incidence however it turns out to be
surprisingly rich.

There is a canonical way of associating a geometry to the pregeometry $\cl$ on
$X$. We define an equivalence relation $\sim$ on $\hat X = X \setminus
\cl(\emptyset)$ by $a \sim b$ iff $\cl(a) = \cl(b)$. Then the operator $\hat
\cl$ defined as $\hat \cl(A/\sim) = \{b/\sim : b \in \cl(A)\}$ is a geometry on
$\hat X$.

Let $\cl$ be a pregeometry on $X$. We distinguish the following types of
pregeometries.
\begin {itemize}
\item The pregeometry is called {\em trivial} if $\cl(A) = \bigcup_{a \in A}
	\cl(a)$. In the associated geometry this means that every set is closed.
\item It is called {\em modular} if $\dim(A \cup B) + \dim(A \cap B) = \dim(A)
	+ \dim(B)$ for all finite dimensional closed $A$ and $B$. The
	pregeometry associated with the linear closure in a vector space is
	modular.
\item The pregeometry is called {\em locally modular} if the above modularity
	equation holds provided $\dim(A \cap B) > 0$. This is equivalent to
	$\cl_a$ being modular for every $a \in X \setminus \cl(\emptyset)$. An
	example of locally modular pregeometry that is not modular is the affine
	closure in a vector space. (The closed sets are translates of linear
	subspaces.)
\item A geometry that is modular but not trivial is called {\em projective}. A
	classical result in projective geometry is that every projective
	geometry of dimension at least $4$ comes from a vector space over a
	division ring.
\end {itemize}

A pregeometry on a structure $M$ is a pregeometry $\cl$ on the base set of $M$
such that every automorphism of $M$ is an automorphism of $\cl$ (i.e. preserves
closed sets). The pregeometry $\cl$ on $M$ is called {\em homogeneous} if
$\dim(M)$ is infinite and for every finite $A \subset M$ and $b, c \in M
\setminus \cl(A)$ there is an automorphism of $M$ fixing $A$ pointwise and
taking $b$ to $c$.

\section {Main Result}

From now on let $G$ be a group and $\cl : \mathcal P(G) \to \mathcal P(G)$ be a
homogeneous pregeometry in the above sense. We call a tuple $\bar a \in G^n$
{\em generic} over $A \subset G$ if $\dim(\bar a/A) = n$. For $A, B \subseteq G$
we say that $B$ is {\em $A$-invariant} if $B$ is fixed setwise by all
automorphism of $G$ that fix $A$ pointwise.

We start with some simple observations.

\begin {proposition}
If $a$ is $A$-invariant for some finite $A \subset G$, then $a \in \cl(A)$.
\end {proposition}

\begin {proof}
Indeed, pick $c \in G$ generic over $A \cup a$. This is possible since $G$ is
infinite dimensional. If $a \not \in \cl(A)$, then there is an automorphism
fixing $A$ and taking $a$ to $c$.
\end {proof}

In particular if $a$ is first-order definable from $A$, then $a \in \cl(A)$. We
will mostly use the Proposition in this form.

\begin {corollary}
Let $A \subset G$ be a finite subset. If $a$ is generic over $A \cup b$, then $a
\cdot b$ is also generic over $A \cup b$.
\end {corollary}

\begin {proof}
\label {generic}
Indeed, $a$ is definable from $b$ and $a \cdot b$. Hence $a \in \cl(A,b, a \cdot
b)$. So in particular $a \cdot b \not \in \cl(A,b)$.
\end {proof}

As another corollary we can exclude one type of pregeometry. 

\begin {corollary}
\label {nontrivial}
The pregeometry $(G, \cl)$ is not trivial.
\end {corollary}

\begin {proof}
Pick a generic pair $a, b$ (i.e. $\dim(a,b) = 2$) and consider $a \cdot b$. By
the above $a \cdot b \not \in \cl(a) \cup \cl(b)$. However $a \cdot b \in
\cl(a,b)$ demonstrating that $\cl$ is not trivial.
\end {proof}

Another consequence is that proper invariant (in particular definable) subgroups
are small.

\begin {proposition}
Let $A$ be a finite subset and $H \le G$ an $A$-invariant subgroup. If $H$
contains an element generic over $A$, then $H = G$.
\end {proposition}

\begin {proof}
Let $b \in G$ be an arbitrary element. Pick $a$ generic over $A \cup b$. Since
$H$ is $A$-invariant and contains a generic element, we have that $a \in H$.
(Recall that there is an automorphism over $A$ mapping one generic element to
another.) By Proposition \ref {generic}, $a \cdot b$ is also generic over $A$.
Hence as above $a \cdot b \in H$. Finally since $H$ is a subgroup, we have that
$b \in H$.
\end {proof}

As a corollary we get the following criterion for a group with homogeneous
pregeometry to be commutative.

\begin {proposition}
If a generic pair commutes in a group with a homogeneous pregeometry, then the
group is commutative.
\end {proposition}

\begin {proof}
Let $(a, b)$ be a generic pair and assume that $a$ and $b$ commute. Consider the
centraliser $C_G(a)$ of $a$. It is definable over $a$ and contains a generic
element $b$. So by the previous proposition, we have $C_G(a) = G$. By
homogeneity the same is true for every generic element. Now given $c \in
\cl(\emptyset)$, by the above it commutes with $a$.  Hence $C_G(c)$ is
$c$-definable and contains a generic element and therefore must be the whole
group $G$.
\end {proof}

We will apply the following improved criterion.

\begin {proposition}
\label {cl-com}
Let $(a, b)$ be a generic pair over some finite $A$. If $b \cdot a \in \cl(A, a
\cdot b)$, then $G$ is commutative.
\end {proposition}

\begin {proof}
Assume that $b \cdot a \in \cl(A, a \cdot b)$. By homogeneity it is true for
every generic pair over $A$. Now $(b, b^{-1} \cdot a)$ is a generic pair over
$A$. And so $b^{-1}ab \in \cl(A,a)$. Now consider the set $\{x \in G: x^{-1}ax =
b^{-1}ab\}$. It is definable over $A, a, b^{-1}ab$ and contains a generic
element $b$. Hence it contains all generic elements over $A, a$. Pick $c$
generic over $A, a, b$. Then $c^{-1}ac = b^{-1}ab$ and therefore $(bc^{-1})a =
a(bc^{-1})$.  Since $(a, bc^{-1})$ is a generic pair, by the above $G$ is
commutative.
\end {proof}

So far we have not used local modularity of the pregeometry. So everyting said
so far also applies to groups with arbitrary homogeneous pregeometry. We now use
the full strength of our assumptions for the main result.

{
\begin {theorem}
A group with locally modular homogeneous pregeometry is commutative.
\end {theorem}
\addtocounter{definition}{-1}
}

\begin {proof}
Assume that $\cl$ is a locally modular homogeneous pregeometry on a group $G$.
Pick a finite $A \not \subseteq \cl(\emptyset)$. Then the localisation $\cl_A$
of $\cl$ at $A$ is modular. (Recall that $\cl_A$ is defined as $\cl_A(X) = \cl(A
\cup X)$.) Consider a generic pair $(a, b)$ over $A$. By the previous
proposition $b \cdot a \in \cl_A(a \cdot b)$ implies that $G$ is commutative.
So assume that $b \cdot a \not \in \cl_A(a \cdot b)$. Pick an element $c \not
\in \cl_A(a, b)$. Similarly we assume that $c \cdot a \not \in \cl_A(a \cdot
c)$. Consider the {\em geometry} of the operator $\cl_A$ on $\cl_A(a, b, c)$.
Since the geometry is modular but not trivial (Corollary \ref {nontrivial}), it
is a projective plane. (This essentially means that any two lines intersect.)

We claim that the lines $(b, c)$, $(ab, ac)$ and $(ba, ca)$ meet in a single
point. Assume the opposite. Then these three lines meet in three different
points $d_1, d_2$ and $d_3$ as shown in the diagram below.

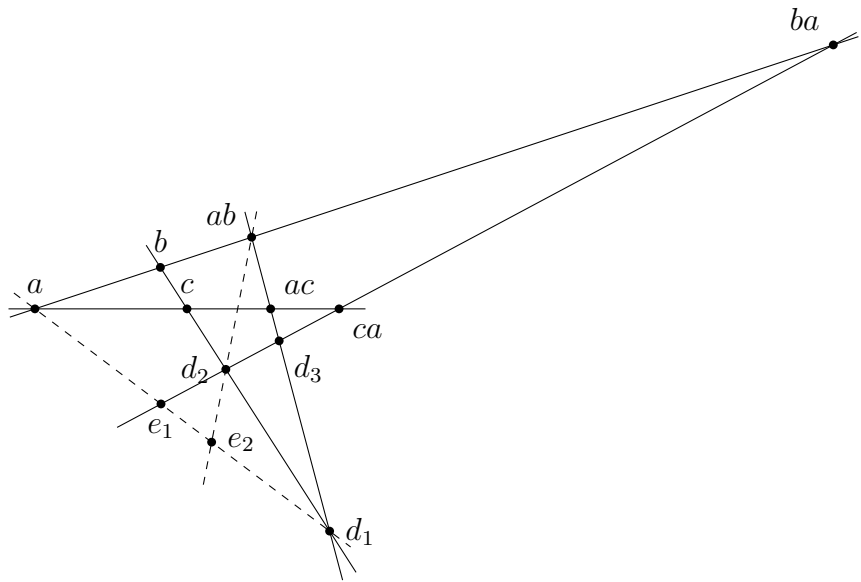
\begin {center}
\begin{tikzpicture}[
	dot/.style 2 args={circle, inner sep=0.1em, fill, label=#2:$#1$, name=#1},
	eline/.style={shorten >=-#1, shorten <=-#1},
	eline/.default={1em},
]

\node [dot={a}{above}] at (0,0) {};
\node [dot={b}{above}] at (1.65,0.55) {};
\node [dot={ab}{160}] at (2.85,0.95) {};
\node [dot={ba}{120}] at (10.5,3.5) {};
\node [dot={c}{above}] at (2,0) {};
\node [dot={ac}{60}] at (3.1,0) {};
\node [dot={ca}{-60}] at (4,0) {};


\draw [eline] (a) -- (ba);
\draw [eline] (a) -- (ca);
\draw [eline, name path=b-to-c] (b) -- ($(c)!-9em!(b)$);
\draw [eline, name path=ab-to-ac] (ab) -- ($(ac)!-8em!(ab)$);
\draw [eline, name path=ba-to-ca] (ba) -- ($(ca)!-7em!(ba)$);

\fill [name intersections={of={b-to-c and ab-to-ac}}] 
	(intersection-1) circle [] node [dot={d_1}{right}] {};
\fill [name intersections={of={b-to-c and ba-to-ca}}] 
	(intersection-1) circle [] node [dot={d_2}{left}] {};
\fill [name intersections={of={ab-to-ac and ba-to-ca}}] 
	(intersection-1) circle [] node [dot={d_3}{-45}] {};

\draw [eline, dashed, name path=D] (a) -- (d_1);

\draw [eline, dashed, name path=ab-to-d2] (ab) -- ($(d_2)!-3em!(ab)$);

\fill [name intersections={of={D and ba-to-ca}}]
	(intersection-1) circle [] node [dot={e_1}{below}] {};
\fill [name intersections={of={D and ab-to-d2}}]
	(intersection-1) circle [] node [dot={e_2}{right}] {};

\end{tikzpicture}
\end {center}

Now consider the line $(a, d_1)$ and assume that it meets the lines $(ba, ca)$
and $(ab, d_2)$ in points $e_1$ and $e_2$ respectively. By the homogeneity
assumption there is an automorphism $f$ of $G$ that fixes $a, e_1, e_2, d_1$ and
takes $b$ to $c$. (Note that $e_1, e_2, d_1$ are points in the geometry so are
equivalence classes of elements of $G$. What we mean by $f$ fixing them is $f$
fixing a representative.) Then $f$ induces a collineation of the projective
plain (which we again denote by $f$). Since $f$ is an automorphism of $G$ it
takes $ab$ to $ac$ and $ba$ to $ca$ respectively. Then $f$ fixes the lines $(b,
c)$ and $(ba, ca)$ (since it also fixes two points $d_1$ and $e_1$ on them).
Thus $f$ fixes their intersection point $d_2$. But this implies that $f$ fixes
the line $(d_2, e_2)$. However $ab$ is on that line whereas $ac$ is not. This
contradiction proves the claim.

But now note that $c^{-1}b \in \cl_A(b, c) \cap \cl_A(ab, ac)$ and $bc^{-1} \in
\cl_A(b, c) \cap \cl_A(ba, ca)$. Thus $c^{-1}b \in \cl_A(bc^{-1})$. Since $(b,
c^{-1})$ is a generic pair over $A$, by Proposition \ref{cl-com} the group $G$
is commutative.
\end {proof}

As a final remark, note that under the assumption of local modularity the
geometry of $\cl$ on $G$ (after localisation) is isomorphic to the projective
geometry of some vector space over a division ring. The group structure of $G$
however need not come from this vector space. Examples of locally modular
strongly minimal groups that are not vector spaces can be found in
\cite{pillay-tori}.

\bibliographystyle{plainnat}
\bibliography{../all}

\begin{thebibliography}{7}
\providecommand{\natexlab}[1]{#1}
\providecommand{\url}[1]{\texttt{#1}}
\expandafter\ifx\csname urlstyle\endcsname\relax
  \providecommand{\doi}[1]{doi: #1}\else
  \providecommand{\doi}{doi: \begingroup \urlstyle{rm}\Url}\fi

\bibitem[Gogacz and Krupi\'nski(2014)]{gogacz}
Tomasz Gogacz and Krzysztof Krupi\'nski.
\newblock On regular groups and fields.
\newblock \emph{The Journal of Symbolic Logic}, 79\penalty0 (3):\penalty0
  826--844, 2014.

\bibitem[Hyttinen(2002)]{hyttinen}
Tapani Hyttinen.
\newblock Groups acting on geometries.
\newblock In Yi~Zhang, editor, \emph{Logic and Algebra}, pages 221--233.
  American Mathematical Society, 2002.

\bibitem[Hyttinen et~al.(2005)Hyttinen, Lessmann, and Shelah]{hls}
Tapani Hyttinen, Olivier Lessmann, and Saharon Shelah.
\newblock Interpreting groups and fields in some nonelementary classes.
\newblock \emph{Journal of Mathematical Logic}, 5\penalty0 (1):\penalty0 1--47,
  2005.

\bibitem[Maesono(2007)]{maesono}
Hisatomo Maesono.
\newblock On quasi-minimal $\omega$-stable groups.
\newblock \emph{RIMS K\^{o}ky\^{u}roku}, 1555:\penalty0 70--72, 2007.

\bibitem[Pillay(1996)]{pillay-tori}
Anand Pillay.
\newblock Definable sets in generic complex tori.
\newblock \emph{Annals of Pure and Applied Logic}, 77:\penalty0 75--80, 1996.

\bibitem[Pillay and Tanovi\'{c}(2011)]{pilovic}
Anand Pillay and Predrag Tanovi\'{c}.
\newblock Generic stability, regularity, and quasiminimality.
\newblock In Bradd Hart, Thomas~G. Kucera, Anand Pillay, and Philip~J. Scott,
  editors, \emph{Models, Logics, and Higher-Dimensional Categories}, pages
  189--212. American Mathematical Society, 2011.

\bibitem[Reineke(1975)]{reineke}
Joachim Reineke.
\newblock Minimale {G}ruppen.
\newblock \emph{Zeitschrift f\"ur Mathematische Logik und Grundlagen der
  Mathematik}, 21\penalty0 (1):\penalty0 357--359, 1975.

\end{thebibliography}
\end {document}